\documentstyle[12pt]{article}

\def\no{\nonumber}
\def\be{\begin{equation}}
\def\ee{\end{equation}}
\def\ba{\begin{eqnarray}}
\def\ea{\end{eqnarray}}
\def\tilde{\widetilde}
\def\btu{\Delta}
\def\btd{\nabla}
\def\e1{\epsilon}
\def\o1{\omega}
\def\01{\Omega}
\def\c1{\gamma}
\def\al{\alpha}

\def\g1{\Sigma}
\def\lmd{\lambda}
\def\l1{\Lambda}
\def\v1{\varphi}
\def\d1{\delta}
\def\part{\partial}
\def\fr{\frac}
\def\f2{F}
\def\h2{{\bf H}}
\def\a2{{\bf A}}
\def\x2{{\bf X}}
\def\t1{\theta}
\def\b1{\beta}
\def\bar{\overline}
\def\bs{\begin{eqnarray*}}
\def\es{\end{eqnarray*}}
\def\p{\partial}
\def\m1{\Theta}
\def\w1{\wedge}
\def\la{\langle}
\def\ra{\rangle}

\begin{document}

\title{Singularities of Codimension Two Mean Curvature Flow of
Symplectic Surfaces}

\author{\begin{tabular}{ccc}
Jingyi Chen\thanks{Chen is supported partially by a Sloan fellowship and a NSERC grant.} &Jiayu
Li\thanks{Li is partially supported by the National Science Foundation of China 
and by the Partner Group of MPI for Mathematics.}\\
Department of Mathematics & Institute of Mathematics\\
The University of British Columbia & Fudan University, Shanghai\\
Vancouver, B.C. V6T 1Z2 & Academia Sinica, Beijing\\
Canada & P.R. China\\
jychen@math.ubc.ca&lijia@math.ac.cn
\end{tabular}}

\date{}

\maketitle


\newtheorem{theorem}{Theorem}[section]
\newtheorem{lemma}[theorem]{Lemma}
\newtheorem{corollary}[theorem]{Corollary}
\newtheorem{remark}[theorem]{Remark}
\newtheorem{definition}[theorem]{Definition}
\newtheorem{proposition}[theorem]{Proposition}
\newtheorem{example}[theorem]{Example}

\section{Introduction}

Geometers have been interested in constructing minimal surfaces
for long time. One possible way to produce such surfaces is to
deform a given surface by its mean curvature vector. More
precisely, the surface evolves in the gradient flow of the area
functional, and such a flow is the so-called mean curvature flow.
A mean curvature flow, however, develops singularity after finite
time in general [H2]. It is therefore desirable to understand
behavior of the flow near the singular points (cf. [B], [CL],
[E1-E2], [H1-H3], [HS1-HS2], [I1-I2], [Wa1], [Wh1-Wh2] and so on).
We consider this problem, in this paper, for compact symplectic
surfaces moving by mean curvature flow in a K\"ahler-Einstein
surface. It was observed in [CT2] that symplectic surface remains
symplectic along the flow (also see [CL], [Wa1]), and the flow has
long time existence in the graphic case as discussed in [CLT] and
[Wa2]. One of our motivations of considering symplectic surfaces
is inspired by the symplectic isotopic problem for symplectic
surfaces in Del Pezzo surfaces, i.e., those complex surfaces with
positive first Chern class. It was conjectured in [T] that every
embedded orientable closed symplectic surface in a compact
K\"ahler-Einstein surface is isotopic to a symplectic minimal
surface in a suitable sense. When K\"ahler-Einstein surfaces are
of positive scalar curvature, this was proved for lower degrees by
using pseudo-holomorphic curves (cf. [ST], [Sh]). It would be
interesting to have a proof of this result, even for lower
degrees, by using the mean curvature flow. In the negative scalar
curvature case, Arezzo [Ar] pointed out that a symplectic minimal
surface may not be holomorphic even it represents a $(1,1)$-type
class, by constructing examples of symplectic minimal surfaces
which are not holomorphic. Moving Lagrange surfaces by mean
curvature flow was studied by Smoczyk [Sm1-Sm3].

Let $M$ be a compact K\"ahler-Einstein surface. On $M$, let
$\omega$ be the K\"ahler form, $\langle,\rangle$ the K\"ahler
metric and $J$ the complex structure which is compatible with
$\omega$ and $\langle,\rangle$, i.e., for any tangent vectors
$X,Y$ to $M$ at an arbitrary point in $M$,
\begin{eqnarray*}
\langle X,Y\rangle&=&\omega(X,JY),\\
\omega(X,Y)&=&\omega(JX,JY).
\end{eqnarray*}

For a compact oriented real surface $\g1$ without boundary which
is smoothly immersed in $M$, one defines, following [CW], the
K\"ahler angle $\alpha$ of $\g1$ in $M$ by
$$
\omega|_\g1 = \cos\alpha \,d\mu_\g1,
$$
where $d\mu_\g1$ is the area element of $\g1$ in the induced
metric from $\langle,\rangle$. As a function on $\g1$, $\alpha$ is
continuous everywhere and is smooth possibly except at the complex
or anti-complex points of $\g1$, i.e. where $\al=0$ or $\pi$. We
say that $\g1$ is a holomorphic curve if $\cos\alpha \equiv 1$,
$\g1$ is a Lagrangian surface if $\cos\alpha \equiv 0$ and $\g1$
is a symplectic surface if $\cos\alpha > 0$.

Given an immersion $F_0:\g1\to M$, we consider a one-parameter
family of immersions $\f2_t=\f2(\cdot ,t):\g1\to M$, and denote
the image surfaces by $\g1_t=\f2_t(\g1)$. The immersed surfaces
$\g1_t$ satisfy a mean curvature flow if
\begin{equation}\label{meaneqn}
\left\{\begin{array}{lll}
\displaystyle \fr{d}{dt}\f2(x,t)=\h2(x,t)\\
\displaystyle \,\,\,\,\,\,\f2(x,0)=\f2_0(x),
\end{array}\right.
\end{equation}
where $\h2(x,t)$ is the mean curvature vector of $\g1_t$ at
$\f2(x,t)$ in $M$.

The standard parabolic theory implies that the mean curvature flow (\ref{meaneqn})
has a smooth solution for short time. More precisely, there exists
$T>0$ such that (\ref{meaneqn}) has a smooth solution in the time
interval $[0, T)$. If the second fundamental form $|\a2|^2$ on
$\g1_t$ is bounded uniformly in $t$ near $T$, then the solution
can be extended smoothly to $[0,T+\e1)$ for some $\e1>0$. However, in
general $\max_{\g1_t} |\a2|^2$ becomes unbounded as $t\to T$. In
this case we say that the mean curvature flow blows up at $T$;
moreover, to classify the singularities of mean curvature flows,
Huisken, according to the blowing up rate of $|\a2|$, introduced

\begin{definition} {\em  (Huisken [H1])
On a mean curvature flow $\g1_t$, suppose that
$$
\lim_{t\to T^-}\max_{\g1_t}|\a2|^2=\infty. $$  If there exists a
positive constant $C$ such that
$$
\limsup_{t\to T^-}\left((T-t)\max_{\g1_t} |\a2|^2\right)\leq C,
$$
the mean curvature flow $\f2$ has a {\it Type I singularity at $T$};
otherwise it has a {\it Type II singularity at $T$}.}
\end{definition}

In the codimension one case, singularities of mean curvature flows
have been studied in depth (cf. [E1], [H1-H3], [HS1-HS2], [I1-I2], [Wh2]).
For higher codimesions, if the initial compact
surface is symplectic in a compact K\"ahler-Einstein surface, the
motion of the mean curvature flow preserves symplecticity of
$\g1_t$ as long as the smooth solution exists; and furthermore the
flow does not develop any Type I singularities (cf. [CT2], [CL], [Wa1]).

In this paper, we shall study the Type II singularities of the
mean curvature flow of a compact symplectic surface in a compact
K\"ahler-Einstein surface. Especially, we shall focus on the
structure of tangent cones of the mean curvature flow where a
singularity occurs at the first singular time $T<\infty$.

To describe the tangent cones, suppose $(X_0,T)$ is a singular
point of the flow (\ref{meaneqn}), i.e. $|\a2(x,t)|$ becomes
unbounded when $(x,t)\to (X_0,T)$. For an arbitrary sequence of
numbers $\lmd\to\infty$ and any $t<0$, if $T+\lambda^{-2}t>0$ we set
$$
F_\lmd(x,t)=\lmd (F(x,T+\lmd^{-2}t)-X_0).
$$
We denote the scaled surface by $(\g1_t^\lmd,d\mu_t^\lmd)$. If the
initial surface is symplectic, it is proved in Lemma
\ref{exit} that there is a subsequence $\lmd_i\to\infty$ such
that for any $t<0$, $(\g1_t^{\lmd_i},d\mu_t^{\lmd_i})$ converges
to $(\g1^{\infty},d\mu^{\infty})$ in the sense of measures;
the limit $\g1^\infty$ is called a {\it tangent cone arising from
the rescaling $\lmd$}, or simply a {\it $\lmd$ tangent cone
at $(X_0,T)$}. This tangent cone is
independent of $t$ as shown in Lemma \ref{exit}. There is also a time
dependent scaling which we would like to consider
\begin{equation}
\tilde{F}(\cdot ,s)=\fr{1}{\sqrt{2(T-t)}}F(\cdot ,t),
\end{equation}
where $s=-\fr{1}{2}\log (T-t)$, $c_0\leq s<\infty $. Here we have chosen the
coordinates so that $X_0=0$. Rescaling of
this type arises naturally in classifying Type I and Type II
singularities for mean curvature flows [H2].  Denote $\tilde{\g1}_s$
the rescaled surface by $\tilde{F}(\cdot,s)$. If a subsequence of
$\tilde{\g1}_s$ converges in measures to a limit $\tilde{\g1}_\infty$,
then the limit is called a
{\it tangent cone arising from the time dependent scaling at
$(X_0,T)$}, or simply a {\it $t$ tangent cone}. In this paper, a
{\it tangent cone}  of the mean curvature flow at $(X_0,T)$ means
either a $\lmd$ tangent cone or a $t$ tangent cone at $(X_0,T)$.

The main result of this paper is

\begin{theorem}\label{main} Let $M$ be a compact K\"ahler-Einstein surface.
If the initial compact surface is symplectic and $T>0$ is the first
blow-up time of the mean curvature flow (\ref{meaneqn}) and $(X_0,T)$
is a singular point, then the tangent cone of the mean curvature flow
at $(X_0,T)$ consists of a finite union of more than one 2-planes
in ${\bf R}^4$ which are complex in a complex structure on ${\bf
R}^4$.
\end{theorem}

It is interesting to compare our result with the one obtained by
Morgan in [M] (also see [MW]): the tangent cone to an oriented
area minimizing surface in ${\bf R}^4$ consists of 2-planes which
are complex in ${\bf R}^4$ (see Remark \ref{morgan}).

An important technical device in the present paper is a
monotonicity formula, which is established in [CL] and also see
Proposition \ref{mono}, for
$$
\int_{\Sigma_t}\frac{1}{\cos\al}\rho(F,t)\phi d\mu_t,
$$
where $F:\Sigma_0\times[0,T)\to M$ is the mean curvature flow,
$\phi$ is a cut-off function supported in a small geodesic ball in
$M$ which is contained in a single coordinate chart, and $\rho$ is
defined by the backward heat kernel on ${\bf R}^4$. The weight
function $1/\cos\alpha$ captures some key geometric information of
the symplectic surfaces in the flow, and this quantity is
sensitive to orientation.

The organization of the paper is follows. From section 2 to
section 5 of this paper, we consider the $\lmd$ scaling near a
singular point $(X_0,T)$. In section 2, first we derive four
integral estimates in Proposition \ref{main2.1}, by using the
weighted monotonicity formula (\ref{mon2}). These estimates play a
crucial role in our analysis of the tangent cones. Then we prove
in Lemma \ref{exit} that a subsequence of the rescaled surfaces
$\g1_t^\lmd$ converges in measure to a limit $\g1^\infty$, which
is independent of $t$. In section 3, we prove that the tangent
cone $\g1^\infty$ is rectifiable. This is done by first showing that the
density function of the tangent cone exists and then applying a
local regularity theorem of White for classical mean curvature
flows (cf. [Wh1], [E1]) to obtain a positive lower bound for the
density function of $\g1^\infty$, and then using Priess's theorem
in [P] to conclude the rectifiability of $\g1^\infty$. The
positive lower bound can also be obtained by using the
isoperimetric inequality (cf. the proof of Lemma \ref{texit}). In
section 4, using Proposition \ref{main2.1} and Allard's
compactness theorem in [A] we first show, in Proposition
\ref{stat}, that $\g1^\infty$ is stationary. Then we prove that
the restriction of the K\"ahler form on the scaled surface
converges to the restriction of a constant 2-form $\omega_0$ in
${\bf R}^4$ along $\g1^\infty$, and further we prove that
$\omega_0|_{\g1^\infty}=\theta_0 d\mu^\infty$ for an ${\cal H}^2$
a.e. constant function $\theta_0$ on $\g1^\infty$. A theorem
of Harvey-Shiffman in [HS] then implies that $\g1^\infty$ is a holomorphic
subvariety. In section 5, we further use Proposition \ref{main2.1}
and Allard's compactness theorem to conclude $F_\infty^\perp$
equals zero and then show $\g1^\infty$ is flat away from its
singular locus, which is a collection of finitely many points.
Section 6 concerns with the time dependent
scaling. The arguments essentially proceed in a similar fashion as those in
the previous sections.

\section{Monotonicity formula and integral estimates}

Let $\g1_t=F(\g1,t)$ be the family of immersed surfaces, which are
determined by the mean curvature flow $F$, in the 4-dimensional
manifold $M$. Denote the Riemannian metric on $M$ by
$\langle\cdot,\cdot\rangle$. In a normal coordinate chart around a
point in $\g1_t$, the induced metric on $\Sigma_t$ from
$\langle\cdot,\cdot\rangle$ is given by
$$
g_{ij}=\la\p_iF,\p_jF\ra,
$$
where $\p_i$ ($i=1,2$) are the partial derivatives with respect to
the local coordinates. In the sequel, we denote by $\btu$ and
$\nabla$ the Laplace operator and covariant derivative for the
induced metric on $\Sigma_t$ respectively. We choose a local field
of orthonormal frames $e_1,e_2$, $v_1,v_2$ of $M$ along
$\Sigma_t$ such that $e_1,e_2$ are tangent vectors of $\g1_t$ and
$v_1,v_2$ are in the normal bundle over $\g1_t$. The second fundemantal
form $\a2$ and the mean curvature vector $\h2$ of $\g1_t$ can be expressed, in the
local frame, as
\begin{eqnarray*}
\a2&=&A^\al v_\al\\
\h2&=&-H^\al v_\al
\end{eqnarray*}
where and throughout this paper all repeated indices are summed over suitable range.
For each $\alpha$, the coefficient $A^\al$ is a $2\times 2$ matrix $(h^\al_{ij})$. By
Weingarten's equation (cf. [Sp]), we have
$$
h_{ij}^\al=\la\p_i v_\al,\p_j F\ra =\la\p_j
v_\al,\p_i\f2\ra=h^\al_{ji}.
$$
The trace and the norm of the second fundamental form of
$\Sigma_t$ in $M$ are:
$$
H^\al = g^{ij}h^\al_{ij}=h_{ii}^\al $$ $$ |\a2|^2=\sum_{\al}
|A^\al |^2=g^{ij}g^{kl}h^\al_{ik}h^\al_{jl}=h_{ik}^\al h_{ik}^\al
.
$$
The area element of the induced metric $g_{ij}$ on $\Sigma_t$ is $\sqrt{\det(g_{ij})}dxdy$. Along
the mean curvature flow, it is well known that
$$
\frac{d}{dt}\sqrt{\det(g_{ij})}=-|\h2|^2 \sqrt{\det(g_{ij})}.
$$
Logarithmic integration implies that ${F}$ remains immersed as
long as the smooth solution of (\ref{meaneqn}) exists.

Let $J_{\g1_t}$ be an almost complex structure in a tubular
neighborhood of $\g1_t$ on $M$ with
\begin{equation}\label{eq2}
\left\{\begin{array}{clcr} J_{\g1_t}e_1&=&e_2\\
J_{\g1_t}e_2&=&-e_1\\ J_{\g1_t}v_1&=&v_2\\
J_{\g1_t}v_2&=&-v_1.
\end{array}\right.
\end{equation}

It is not difficult to verify (cf. Lemma 3.1 in [CL]), with $\overline{\nabla}$ being the
covariant derivative of the metric $\langle\cdot,\cdot\rangle$ on $M$, that
\begin{eqnarray}\label{djh}
|\overline{\nabla}J_{\g1_t}|^2&=&
|h_{11}^2+h_{12}^1|^2+|h^2_{21}+h^1_{22}|^2
+|h_{12}^2-h_{11}^1|^2+|h^2_{22}-h^1_{21}|^2 \nonumber\\
&=&\frac{1}{2}|\h2|^2+\frac{1}{2}\left(((h^1_{11}+h^1_{22})+2(h^2_{12}-h^1_{22}))^2
+(h^2_{11}+h^2_{22}+2h^1_{21}-2h^2_{11})^2\right)\nonumber\\
&\geq&
\fr{1}{2}|\h2|^2.
\end{eqnarray}

Let $H (\x2,\x2_0 ,t)$ be the backward heat kernel on ${\bf R}^4$.
Define
$$
\rho (\x2,\x2_0,t,t_0)=4\pi (t_0-t)H (\x2,\x2_0 ,t)=\fr{1}{4\pi
(t_0-t)}\exp \left(-\fr{|\x2-\x2_0|^2} {4(t_0-t)}\right)
$$
for $t<t_0$. Let $i_M$ be the injective radius of $M^4$. We choose a cut off
function $\phi\in C_0^\infty(B_{2r}(\x2_0))$ with $\phi\equiv 1$
in $B_{r}(\x2_0)$, where $\x2_0\in M$, $0<2r<i_M$. Choose a normal
coordinates in $B_{2r}(\x2_0)$ and express $\f2$, by the
coordinates $(F^1,F^2,F^3,F^4)$, as a surface in ${\bf R}^4$. We
define
\begin{equation}\label{Phi}
\Phi (\x2_0,t_0,t)=\int_{\g1_t}\phi (\f2 )\rho
(\f2,\x2_0,t,t_0)d\mu_t.
\end{equation}

Huisken derived the following useful monotonicity formula in [H1]:
there are positive constants $c_1$ and $c_2$ depending only on
$M^4$, $\f2_0$ and $r$ where $r$ is the constant in the definition
of $\phi$, such that \begin{eqnarray}\label{mon1} \fr{\p}{\p
t}\left(e^{c_1\sqrt{t_0-t}}\Phi (\x2_0,t_0,t)\right)&\leq&
-e^{c_1\sqrt{t_0-t}}\int_{\g1_t}\phi\rho (\f2,\x2_0
,t,t_0)\left|\h2+\fr{(\f2-\x2_0)^{\perp}}
{2(t_0-t)}\right|^2d\mu_t \nonumber\\
&&+c_2(t_0-t). \end{eqnarray} Note that $c_1$ and $c_2$ equal 0
when $M$ is a Euclidean space.

When $M$ is a K\"ahler-Einstein surface with scalar curvature $R$
and $\Sigma_t$ evolves under the mean curvature flow, the K\"ahler
angle $\alpha$ of $\g1_t$ in $M$ satisfies the parabolic equation (cf. [CL], [CT]):
\be\label{kef} \left(\fr{\p}{\p t}-\btu \right)\cos\al =
|\overline{\nabla}J_{\g1_t}|^2\cos\al + R\sin^2\al\cos\al. \ee

Suppose that the initial surface is symplectic, i.e., $\cos
\al(\cdot,0)$ has a positive lower bound. Then by applying the
parabolic maximum principle to the evolution equaiotn (\ref{kef}),
one concludes that $\cos\alpha$ remains positive
as long as the mean curvature flow has a smooth solution, no matter $R$ is
positive, 0 or negative (cf. [CT2], [CL], [Wa1]).

Let $R_0=\max\{0,-R\}$ and set
$$
v(x,t)=e^{R_0t}\cos\al(x,t).
$$
By (\ref{kef}), we have
\begin{equation}
\left(\fr{\p}{\p t}-\btu \right)\fr{1}{v}\leq
-|\overline{\nabla}J_{\g1_t}|^2\fr{1}{v}-\fr{2}{v^3}|\btd v|^2.
\end{equation}
Along the flow, we introduce a function
\begin{equation}\label{Psi}
\Psi (\x2_0,t_0,t)=\int_{\g1_t}\fr{1}{v}\phi\rho
(\f2,\x2_0,t,t_0)d\mu_t.
\end{equation}

The following weighted monotonicity formula in [CL] will play a
crucial role in this paper.
\begin{proposition} \label{mono}(Weighted Monotonicity Formula) If the initial compact
surface $\Sigma_0$ is symplectic in a K\"ahler-Einstein surface $M$ and $\Sigma_t$
evolves under the mean curvature flow (\ref{meaneqn}), then
\ba \label{mon2}
\fr{\p}{\p t}\left(e^{c_1\sqrt{t_0-t}}\Psi (X_0,t_0,t)\right)
&\leq&-e^{c_1\sqrt{t_0-t}} \left( \int_{\g1_t}\fr{1}{v}\phi\rho
(\f2,\x2_0,t,t_0)
\left|\h2+\fr{(\f2-\x2_0)^{\perp}}{2(t_0-t)}\right|^2d\mu_t \right. \no\\
&&\left.+ \int_{\g1_t}\fr{1}{v}\phi\rho
(\f2,\x2_0,t,t_0)\left|\overline{\nabla}J_{\g1_t}\right|^2d\mu_t\right.\nonumber\\
&&\left.+\int_{\g1_t}\fr{2}{v^3}\left|\btd v\right|^2\phi\rho
(\f2,\x2_0,t,t_0)d\mu_t \right)+c_2(t_0-t). \ea Here the positive
constants $c_1$ and $c_2$ depend on $M$, $\f2_0$ and $r$ where $r$
is the constant in the definition of $\phi$.
\end{proposition}

Suppose that $(X_0,T)$ is a singular point of the mean curvature
flow (\ref{meaneqn}). We now describe the rescaling process around
$(X_0,T)$. For any $t<0$, we set
$$
F_\lmd(x,t)=\lmd (F(x,T+\lmd^{-2}t)-X_0),
$$
where $\lmd$ are positive constants which go to infinity. The
scaled surface is denoted by $\g1_t^\lmd=F_\lmd(\Sigma,t)$ on
which $d\mu^\lmd_t$ is the area element obtained from $d\mu_t$. If
$g^\lmd$ is the metric on $\Sigma^\lmd_t$, it is clear that
$$
g_{ij}^\lmd=\lmd^2g_{ij},~~~~(g^\lmd)^{ij}=\lmd^{-2}g^{ij}.
$$
We therefore have
\begin{eqnarray}
\fr{\p\f2_\lmd}{\p t}&=&\lmd^{-1}\fr{\p\f2}{\p t}\\
\h2_\lmd &=&\lmd^{-1}\h2\\
|\a2_\lmd|^2&=&\lmd^{-2}|\a2|^2.\end{eqnarray} It follows that the
scaled surface also evolves by a mean curvature flow
\begin{equation}
\fr{\p\f2_\lmd}{\p t}=\h2_\lmd .
\end{equation}

\begin{proposition}
\label{main2.1} Let $M$ be a K\"ahler-Einstein surface. If the
initial compact surface is symplectic, then for any $R>0$ and any
$-\infty<s_1<s_2<0$, we have \be\label{m2.1}
\int_{s_1}^{s_2}\int_{\g1_t^\lmd\cap B_R(0)}
|\bar{\btd}J_{\g1_t^\lmd}|^2d\mu_t^\lmd dt\to 0~~{\rm
as}~~\lmd\to\infty, \ee \be\label{m2.2}
\int_{s_1}^{s_2}\int_{\g1_t^\lmd\cap B_R(0)} |\btd \cos
\al_{\lmd}|^2d\mu_t^\lmd dt\to 0~~{\rm as}~~\lmd\to\infty, \ee
\be\label{m2.3} \int_{s_1}^{s_2}\int_{\g1_t^\lmd\cap
B_R(0)}|\h2_\lmd|^2d\mu_t^\lmd dt\to 0~~{\rm as}~~\lmd\to\infty,
\ee and \be\label{m2.4} \int_{s_1}^{s_2}\int_{\g1_t^\lmd\cap
B_R(0)}|F_\lmd^{\perp}|^2d\mu_t^\lmd dt\to 0~~{\rm
as}~~\lmd\to\infty. \ee
\end{proposition}
{\it Proof:} For any $R>0$, we choose a cut-off function
$\phi_R\in C_0^{\infty}(B_{2R}(0))$ with $\phi_R\equiv 1$ in
$B_R(0)$, where $B_r(0)$ is the metric ball centered at $0$ with
radius $r$ in ${\bf R}^4$. For any fixed $t<0$, the mean curvature flow (\ref{meaneqn}) has a
smooth solution near $T+\lmd^{-2}t<T$ for sufficiently large $\lmd$, since
$T>0$ is the first blow-up time of the flow.
It is clear \bs
\lefteqn{\int_{\g1_t^\lmd}\fr{1}{v_\lmd}\fr{1}{0-t}\phi_R(\f2_\lmd
)\exp
\left(-\fr{|\f2_\lmd |^2}{4(0-t)}\right)d\mu_t^\lmd}\\
&=&\int_{\g1_{T+\lmd^{-2}t}}\fr{1}{v_\lmd}\phi (\f2_\lmd)
\fr{1}{T-(T+\lmd^{-2}t)}\exp
\left(-\fr{|\f2(x,T+\lmd^{-2}t)-X_0|^2}{4(T-(T+\lmd^{-2}t))}\right)d\mu_t
, \es where $\phi$ is the function defined in the definition of
$\Phi$. Note that $T+\lmd^{-2}t\to T$ for any fixed $t$ as $\lmd\to\infty$. By
(\ref{mon2}),
$$
\frac{\partial}{\partial t}\left(e^{c_1\sqrt{t_0-t}}\Psi\right)
\leq c_2(t_0-t),
$$
and it then follows that $\lim_{t\to t_0}e^{c_1\sqrt{t_0-t}}\Psi$
exists. This implies, by taking $t_0=T$ and $t=T+\lmd^{-2}s$, that
for any fixed $s_1$ and $s_2$ with $-\infty <s_1<s_2<0$,
\begin{eqnarray}\label{limit} \lefteqn{e^{c_1\sqrt{T-(T+\lmd^{-2}s_2)}}
\int_{\g1_{s_2}^\lmd}\fr{1}{v_\lmd}\phi_R\fr{1}{0-s_2}\exp
\left(-\fr{|\f2_\lmd|^2}{4(0-s_2)}\right)d\mu_{s_2}^\lmd}\nonumber\\
&&-e^{c_1\sqrt{T-(T+\lmd^{-2}s_1)}}\int_{\g1_{s_1}^\lmd}\fr{1}{v_\lmd}
\phi_R\fr{1}{0-s_1}\exp
\left(-\fr{|\f2_\lmd |^2}{4(0-s_1)}\right)d\mu_{s_1}^\lmd \nonumber\\
& \to &0 ~~\hbox{as}~~ \lmd\to\infty. \end{eqnarray} Integrating
(\ref{mon2}) from $s_1$ to $s_2$ yields \ba\label{djh2}
\lefteqn{-e^{c_1\sqrt{-\lmd^{-2}s_2}}\int_{\g1_{s_2}^\lmd}\fr{1}{v_\lmd}
\phi_R\fr{1}{0-s_2}\exp
\left(-\fr{|\f2_\lmd|^2}{4(0-s_2)}\right)d\mu_{s_2}^\lmd}\no\\
&&+e^{c_1\sqrt{-\lmd^{-2}s_1}}\int_{\g1_{s_1}^\lmd}\fr{1}{v_\lmd}\phi_R\fr{1}{0-s_1}\exp
\left(-\fr{|\f2_\lmd|^2}{4(0-s_1)}\right)d\mu_{s_1}^\lmd\no \\
&\geq &
\int_{s_1}^{s_2}e^{c_1\sqrt{-\lmd^{-2}t}}\int_{\g1_t^\lmd}\fr{1}{v_\lmd}\phi_R\rho
(\f2_\lmd ,t)\left|\h2_k+\fr{(\f2_k)^{\perp}}{2(t_0-t)}\right|^2d\mu_t^\lmd\no\\
&&+\int_{s_1}^{s_2}e^{c_1\sqrt{-\lmd^{-2}t}}\int_{\g1_t^\lmd}\fr{1}{v_\lmd}\phi_R\rho
(\f2_k,t)|\overline{\nabla}J_{\g1_t^\lmd}|^2d\mu_t^\lmd\no\\
&&+\int_{s_1}^{s_2}e^{c_1\sqrt{-\lmd^{-2}t}}\int_{\g1_t^\lmd}\fr{2}{v_\lmd^3}|\btd
v_\lmd|^2\phi_R\rho (\f2_\lmd
,t)d\mu_t^\lmd\no\\
&&-c_2\lmd^{-2}(s_2^2-s_1^2). \ea
Putting (\ref{limit}) and (\ref{djh2}) together, we have
$$
\lim_{\lmd\to\infty}\int_{s_1}^{s_2}\int_{\g1_t^\lmd}\phi_R\rho
(\f2_k,t)\left|\overline{\nabla}J_{\g1_t^\lmd}\right|^2d\mu_t^\lmd=0,
$$
and
$$
\lim_{\lmd\to\infty}\int_{s_1}^{s_2}\int_{\g1_t^\lmd}\left|\btd
v_\lmd\right|^2\phi_R\rho (\f2_\lmd ,t)d\mu_t^\lmd= 0,
$$
which yield (\ref{m2.1}) and (\ref{m2.2}) respectively, and
\be\label{m2.5}
\lim_{\lmd\to\infty}\int_{s_1}^{s_2}\int_{\g1_t^\lmd}\phi_R\rho
(\f2_\lmd
,t)\left|\h2_\lmd+\fr{(\f2_\lmd)^{\perp}}{2(t_0-t)}\right|^2d\mu_t^\lmd=
0. \ee
Finally, (\ref{djh}) and (\ref{m2.1}) imply (\ref{m2.3}),
and (\ref{m2.3}) and (\ref{m2.5}) imply (\ref{m2.4}).  \hfill
Q.E.D.

\begin{lemma}\label{exit}
For any $\lmd,R>0$ and any $t<0$, \be\label{fm}
\mu_t^{\lmd}(\g1_t^\lmd\cap B_R(0))\leq CR^2,
 \ee
where $B_R(0)$ is a metric ball in ${\bf R}^4$ and $C>0$ is
independent of $\lmd$. For any sequence $\lmd_i\to\infty$, there
is a subsequence $\lmd_k\to\infty$ such that
$(\g1_t^{\lmd_k},\mu_t^{\lmd_k})\to (\g1^\infty,\mu^\infty)$ in
the sense of measure, for any fixed $t<0$, where
$(\g1^\infty,\mu^\infty)$ is independent of $t$. The multiplicity
of $\g1^\infty $ is finite.
\end{lemma}
{\it Proof:}  We shall first prove the inequality (\ref{fm}). We shall
use $C$ below for uniform positive constants which are independent
of $R$ and $\lmd$. Straightforward computation shows
 \bs
\mu_t^{\lmd}(\g1_t^\lmd\cap B_R(0))&=
&\lmd^2\int_{\g1_{T+\lmd^{-2}t}\cap B_{\lmd^{-1}R}(X_0)}d\mu_t\\
&=& R^2(\lmd^{-1}R)^{-2}\int_{\g1_{T+\lmd^{-2}t}\cap B_{\lmd^{-1}R}(X_0)}d\mu_t\\
&\leq & CR^2\int_{\g1_{T+\lmd^{-2}t}\cap B_{\lmd^{-1}R}(X_0)}
\frac{1}{4\pi (\lmd^{-1}R)^2}e^{-\frac{|X-X_0|^2}{4(\lmd^{-1}R)^2}}d\mu_t\\
&=&CR^2\Phi(X_0,T+(\lmd^{-1}R)^2+\lmd^{-2}t, T+\lmd^{-2}t). \es By
the monotonicity inequality (\ref{mon1}), we have \bs
\mu_t^{\lmd}(\g1_t^\lmd\cap B_R(0))&\leq &CR^2\left(\Phi(X_0,
T+(\lmd^{-1}R)^2+\lmd^{-2}t, T/2)+C\right)\\
&\leq &C\fr{R^2}{T}\left(\mu_{T/2}(\g1_{T/2})+C\right). \es Since
$$
\fr{\p}{\p t}\mu_t(\g1_t)=-\int_{\g1_t}|\h2|^2d\mu_t,
$$
we can now conclude (\ref{fm})
$$
\mu_t^{\lmd}(\g1_t^\lmd\cap B_R(0))\leq CR^2.
$$
By (\ref{fm}), the compactness theorem of the measures (c.f.
[Si1], 4.4) and a diagonal subsequence argument, we conclude
that there is a subsequence $\lmd_k\to\infty$ such that
$(\g1_{t_0}^{\lmd_k},\mu_{t_0}^{\lmd_k})\to
(\g1_{t_0}^\infty,\mu_{t_0}^\infty)$ in the sense of measures for a
fixed $t_0<0$.

We now show that, for any $t<0$, the subsequence $\lmd_k$ which we have
chosen above satisfies $(\g1_{t}^{\lmd_k},\mu_{t}^{\lmd_k})\to
(\g1_{t_0}^\infty,\mu_{t_0}^\infty)$ in the sense of measure. And
consequently the limiting surface
$(\g1^\infty_{t_0},\mu^\infty_{t_0})$ is independent of $t_0$.
Recall that the following standard formula for mean curvature flows
\begin{equation}\label{variation}
\frac{d}{dt}\int_{\Sigma^\lmd_t}\phi d\mu^\lmd_t
=-\int_{\Sigma^\lmd_t}\left(\phi|\h2_\lmd
|^2+\nabla\phi\cdot\h2_\lmd\right)d\mu^\lmd_t
\end{equation}
is valid for any test function $\phi\in C_0^\infty(M)$ (cf. (1) in
Section 6 in [I2] and [B] in the varifold setting).

Then for any given $t<0$ integrating (\ref{variation}) yields
 \bs
 \int_{\Sigma^{\lmd_k}_t}\phi d\mu^{\lmd_k}_t
 -\int_{\Sigma^{\lmd_k}_{t_0}}\phi d\mu^{\lmd_k}_{t_0} &=&
 \int_t^{t_0}\int_{\Sigma^{\lmd_k}_t}\left(\phi|\h2_{\lmd_k}
|^2+\nabla\phi\cdot\h2_{\lmd_k}\right)d\mu^{\lmd_k}_tdt\\
&\to & 0~~{\rm as}~~k\to\infty~~{\rm by}~~(\ref{m2.3}) .\es So,
for any fixed $t<0$, $(\g1_{t}^{\lmd_k},\mu_{t}^{\lmd_k})\to
(\g1_{t_0}^\infty,\mu_{t_0}^\infty)$ in the sense of measures as
$k\to\infty$. We denote $(\g1_{t_0}^\infty,\mu_{t_0}^\infty)$ by
$(\g1^\infty,\mu^\infty)$, which is independent of $t_0$.

The inequality (\ref{fm}) yields a uniform upper bound on
$R^{-2}\mu^{\lmd_k}_t(\g1^{\lmd_k}_t\cap B_R(0))$, which
yields finiteness of the multiplicity of $\g1^\infty$.
\hfill Q.E.D.

\begin{definition}\label{bubble1}\label{bb}
{\em Let $(X_0,T)$ be a singular point of the mean curvature flow
of a closed symplectic surface $\Sigma_0$ in a compact
K\"ahler-Einstein surface $M$. We call $(\g1^\infty, d\mu^\infty)$
obtained in Lemma \ref{exit} {\it a
$\lmd$ tangent cone of the mean curvature flow $\g1_t$ at
$(X_0,T)$.}}
\end{definition}

\section{Rectifiability of the $\lmd$ tangent cones}

In this section we shall show that the $\lmd$ tangent cone $\g1^\infty$ is ${\mathcal
H}^2$-rectifiable, where ${\mathcal H}^2$ is the 2-dimensional Hausdorff measure.

\begin{proposition}\label{rectif}
Let $M$ be a compact K\"ahler-Einstein surface. If the initial compact
surface $\g1_0$ is symplectic, then the $\lmd$ tangent cone $(\g1^\infty
,d\mu^\infty )$ of the mean curvature flow at $(X_0,T)$ is
${\mathcal H}^2$-rectifiable.
\end{proposition}
{\it Proof:} We set
$$
A_R=\left\{t\in (-\infty,0)\left| ~\lim_{k\to
\infty}\int_{\g1_t^k\cap B_R(0)}|{\bf H}_k|^2d\mu_t^k
\not=0\right\}, \right.
$$
and
$$
A=\bigcup_{R>0}A_R.
$$

Denote the measures of $A_R$ and $A$ by $|A_R|$ and $|A|$
respectively. It is clear from (\ref{m2.3}) that $|A_R|=0$ for any $R>0$.
So $|A|=0$.

For any $\xi\in\g1^\infty$, choose $\xi_k\in\g1_t^k$
with $\xi_k\to\xi$ as $k\to\infty$. By the monotonicity identity
(17.4) in [Si1], we have \ba\label{smon}
\sigma^{-2}\mu_t^k(B_\sigma(\xi_k))&=&\rho^{-2}\mu_t^k(B_\rho(\xi_k))
-\int_{B_\rho(\xi_k)\setminus B_\sigma(\xi_k)}\fr{|D^\perp r|^2}{r^2}d\mu_t^k\no\\
&& -\fr{1}{2}\int_{B_\rho(\xi_k)}(x-\xi_k)\cdot {\bf H}_k
\left(\fr{1}{r_\sigma^2}-\fr{1}{\rho^2}\right)d\mu_t^k, \ea for
all $0<\sigma \leq \rho$, where $\mu_t^k(B_\sigma(\xi_k))$ is the
area of $\Sigma^k_t\cap B_\sigma((\xi_k))$, $r=r(x)$ is the
distance from $\xi_k$ to $x$, $r_\sigma =\max \{r,\sigma\}$, and
$D^\perp r$ denotes the orthogonal projection of $Dr$ (which is a
vector of length 1) onto $(T_{\xi_k}\g1^k_t)^\perp$. Choosing
$t\not\in A$, we have
$$
\lim_{k\to\infty}\int_{B_\rho(\xi_k)}|{\bf H}_k|^2d\mu_t^k=0.
$$
H\"older's inequality and (\ref{fm}) then lead to
\begin{eqnarray}\label{estimate}
\lefteqn{\lim_{k\to\infty}\left|\int_{B_\rho(\xi_k)}(x-\xi_k)\cdot {\bf H}_k
\left(\fr{1}{r_\sigma^2}-\fr{1}{\rho^2}\right)d\mu_t^k\right|}\nonumber\\
&\leq&
C\rho\left(\fr{1}{\sigma^2}-\fr{1}{\rho^2}\right)\lim_{k\to\infty}
\left(\sqrt{\mu^k_t(B_{\rho}(\xi_k))}
\sqrt{\int_{B_\rho(\xi_k)}|\h2_k|^2 d\mu^k_t}\right)\nonumber\\
&\leq& C\rho^2\left(\fr{1}{\sigma^2}-\fr{1}{\rho^2}\right)
\lim_{k\to\infty}\sqrt{\int_{B_\rho(\xi_k)}|\h2_k|^2 d\mu^k_t}\nonumber\\
&=&0.
\end{eqnarray}
Letting $k\to\infty$ in (\ref{smon}) and using (\ref{estimate}), we obtain
$$
\sigma^{-2}\mu^\infty(B_\sigma(\xi))\leq
\rho^{-2}\mu^\infty(B_\rho(\xi)),
$$
for all $0<\sigma \leq \rho$. By (\ref{fm}) we know that
$$
 \lim_{\rho\to 0}\rho^{-2}\mu^\infty(B_\rho(\xi))<C<\infty.
$$
Therefore, $\lim_{\rho\to 0}\rho^{-2}\mu^\infty(B_\rho(\xi))$
exists.

We shall show that there exists a positive number $r_0$ such that
for any $0<r<r_0<1$ the following density estimate holds
\be\label{rec} \lim_{\rho\to
0}\rho^{-2}\mu^\infty(B_\rho(\xi))\geq \fr{1}{16}. \ee

Assume (\ref{rec}) fails to hold. Then there is $\rho_0>0$ such
that
$$
(2\rho_0)^{-2}\mu^\infty(B_{2\rho_0}(\xi))<\frac{1}{16}.
$$
By the monotonicity formula (\ref{smon}) and that $\mu^k_t$
converges to $\mu^\infty$ as measures, there exists $k_0>0$ such
that, for all $0<\rho<2\rho_0$ and $k>k_0$, we have
\be\label{rec1} \rho^{-2}\mu_t^k(B_\rho(\xi))< \fr{1}{8}.\ee
Take a cut-off function $\phi_\rho\in C^{\infty}_0(B_\rho(\xi))$
on the 4-dimensional ball $B_\rho(\xi_k)$ so that
\begin{eqnarray*}
&&\phi_\rho\equiv 1~~{\rm in}~~B_{\fr{\rho}{2}}(\xi)\\
&&0\leq\phi_\rho\leq 1,~~{\rm and}~~| \btd\phi_\rho
|\leq\fr{C}{\rho},~~{\rm in}~ B_\rho(\xi).
\end{eqnarray*}
From (\ref{variation}), we have\bs
\lefteqn{\rho^{-2}\int_{B_\rho(\xi)}\phi_\rho d\mu_{t-r^2}^k
-\rho^{-2}\int_{B_\rho(\xi)}\phi_\rho d\mu_{t}^k
}\\
&\leq & 2\rho^{-2}\int_{t-r^2}^t\int_{B_{\rho}(\xi)}|{\bf H}_k|
^2d\mu_s^kds+C\rho^{-3}\int_{t-r^2}^t\int_{B_\rho(\xi)}|{\bf H}_k|
d\mu_s^kds\\
&\leq & C\rho^{-2}\int_{t-r^2}^t\int_{B_{\rho}(\xi)}|{\bf H}_k|
^2d\mu_s^kds +
C\rho^{-3}\int_{t-r^2}^t\left(\int_{B_{\rho}(\xi)}|{\bf H}_k|
^2d\mu_s^k\right)^{1/2}\mu_s^k(B_\rho(\xi))^{1/2}ds \\
&\leq &
C\rho^{-2}\int_{t-r^2}^t\int_{B_{\rho}(\xi)}|{\bf H}_k|
^2d\mu_s^kds +
C\rho^{-2}\int_{t-r^2}^t\left(\int_{B_{\rho}(\xi)}|{\bf H}_k|
^2d\mu_s^k\right)^{1/2}ds~~{\rm by}~(\ref{fm})\\
&\to& 0, ~~{\rm as}~~ k\to\infty~~{\rm by}~(\ref{m2.3}). \es Here
we have used $C$ for uniform positive constants which are
independent of $k$ and $\rho$.
Therefore, there are constants $\d1_1>0$ and $k_1>0$ such that for
all $\rho$ and $k$ with $0<\rho<\d1_1$, $0<r<1$, and $k>k_1$ the
estimate \be \label{rec2}
\rho^{-2}\mu_{t-r^2}^k(B_\rho(\xi))<\frac{1}{4}< 1 \ee holds. Let
$d\sigma^k_{t-r^2}$ be the arc-length element of $\partial B_{\rho}(\xi)\cap \g1^k_{t-r^2}$.
By the co-area formula, for $0<r<<1$,
\begin{eqnarray}\label{density}
\Phi_k(\xi,t,t-r^2) &=&\frac{1}{4\pi r^2}\int_{\Sigma_{t-r^2}^k}
e^{-\frac{|F_k-\xi|^2}{4r^2}}d\mu^k_{t-r^2}\nonumber\\
&\leq &\fr{1}{4\pi r^2}\int_0^{\d1_1}\int_{\part
 B_\rho(\xi)\cap\g1_{t-r^2}^k}e^{-\fr{\rho^2}{4r^2}}d\sigma_{t-r^2}^k d\rho\nonumber\\
&&+\fr{1}{4\pi r^2}\int_{\d1_1}^{\infty}\int_{\part
 B_\rho(\xi)\cap\g1_{t-r^2}^k}e^{-\fr{\rho^2}{4r^2}}d\sigma_{t-r^2}^k d\rho\nonumber\\
&\leq &\fr{1}{4\pi r^2}\int_0^{\d1_1}
e^{-\fr{\rho^2}{4r^2}}\frac{d}{d\rho}{\rm
Vol}(B_\rho(\xi)\cap\Sigma^k_{t-r^2})d\rho\nonumber\\
&&+\fr{1}{4\pi r^2}\int_{\d1_1}^\infty
e^{-\fr{\rho^2}{4r^2}}\frac{d}{d\rho}{\rm
Vol}(B_\rho(\xi)\cap\Sigma^k_{t-r^2})d\rho \nonumber\\
&\leq &\fr{1}{4\pi r^2}\int_0^{\d1_1}
e^{-\fr{\rho^2}{4r^2}}\fr{\rho}{2r^2}{\rm Vol}(
 B_\rho(\xi)\cap\g1_{t-r^2}^k)d\rho\nonumber\\
&&+\fr{1}{4\pi r^2}\int_{\d1_1}^\infty
e^{-\fr{\rho^2}{4r^2}}\fr{\rho}{2r^2}{\rm Vol}(
 B_\rho(\xi)\cap\g1_{t-r^2}^k)d\rho\nonumber\\
&\leq &\fr{1}{4\pi r^2}\int_0^{\d1_1}
e^{-\fr{\rho^2}{4r^2}}\fr{\rho^3}{2r^2}d\rho + o(r)~~{\rm
by}~~(\ref{rec2})~~{\rm and}~~(\ref{fm})\nonumber\\
&\leq & 1 +o(r).
\end{eqnarray}

For any classical mean curvature flow $\Gamma_t$ in a compact
Riemannian manifold which is isometrically embedded in ${\bf
R}^N$, White proves a local regularity theorem (Theorem 3.1 and
Theorem 4.1 in [Wh1]). When dim$\Gamma_t=2$, White's theorem
asserts that there is a constant $\epsilon>0$ such that if the
Gaussian density satisfies
$$
\lim_{r\rightarrow 0}\int_{\Gamma_{t-r^2}}\frac{1}{4\pi r^2}
\exp\left(-\fr{|y-x|^2}{4r^2}\right)d\mu(y) <1+\epsilon,
$$
then the mean curvature flow is smooth in a neighborhood of $x$.
Combining this regularity result with (\ref{density}), we are led
to choose $r>0$ sufficiently small and then conclude that
$$
\sup_{B_r(\xi)\cap\g1_t^k}|\a2_k|\leq C
$$
and consequently $\g1_t^k$ converges strongly in $B_r(\xi)\cap\g1_t^k$ to
$\Sigma^\infty_t\cap B_r(\xi)$, as $k\to\infty$. So $\Sigma^\infty\cap B_r(\xi)$ is
smooth. Smoothness of $\Sigma^\infty\cap B_r(\xi)$ immediately implies
$$
\lim_{\rho\to 0}\rho^{-2}\mu^\infty(B_\rho(\xi))= 1.
$$
This contradicts (\ref{rec1}). Hence we have established
(\ref{rec}).

In summary, we have shown that $\lim_{\rho\to
0}\rho^{-2}\mu^\infty(B_\rho(\xi))$ exists and for ${\mathcal
H}^2$ almost all $\xi\in\g1^\infty$,
\begin{equation}\label{den}
\fr{1}{16}\leq\lim_{\rho\to
0}\rho^{-2}\mu^\infty(B_\rho(\xi))<\infty.
\end{equation}

Finally, we recall a fundamental theorem of Priess in [P]: if
$0\leq m\leq n$ are integers and $\Omega$ is a Borel measure on
${\bf R}^n$ such that
$$
0<\lim_{r\rightarrow 0}\fr{\Omega(B_r(x))}{r^m}<\infty,
$$
for almost all $x\in\Omega$, then $\Omega$ is $m$-rectifiable. Now
we conclude from (\ref{den}) that $(\g1^\infty, \mu^\infty)$ is
${\mathcal H}^2$-rectifiable. \hfill
Q.E.D.

\section{Holomorphicity of the $\lmd$ tangent cones}
In this section, we shall first show that the $\lmd$ tangent cone
$\g1^\infty$ is stationary and then prove that $\g1^\infty$ is a complex
subvariety in ${\bf R}^4$. This result allows us to assert that
the set of singular points of $\g1^\infty$ consists of discrete
points.

A $k$-varifold is a Radon measure on $G^k(M)$,
where $G^k(M)$ is the Grassmann bundle of all $k$-planes tangent to $M$.
Allard's compactness theorem for rectifiable varifolds (6.4 in [A], also see
1.9 in [I2] and Theorem 42.7 in [Si1]) can be stated as follows.
\begin{theorem}\label{cmpt}
(Allard's compactness theorem) Let $(V_i,\mu_i)$ be a sequence of rectifiable
$k$-varifolds in $M$ with
$$
\sup_{i\geq 1}(\mu_i(U)+|\d1V_i|(U))<\infty~~{\rm
for~each}~U\subset\subset M.
$$
Then there is a varifold $(V,\mu)$ of locally bounded first variation
and a subsequence, which we
also denote by $(V_i,\mu_i)$, such that

(i)  Convergence of measures: $\mu_i\to\mu$ as Radon
measures on $M$,

(ii) Convergence of tangent planes: $V_i\to V$ as
Radon measures on $G^k(M)$,

(iii) Convergence of first variations: $\d1 V_i\to \d1 V$ as $TM$-valued
Radon measures,

(iv) Lower semicontinuity of total first variations: $|\d1
V|\leq\liminf_{i\to\infty} |\d1 V_i|$ as Radon measures.
\end{theorem}

We first show that the $\lmd$ tangent cone is stationary.

\begin{proposition}\label{stat}
Let $M$ be a compact K\"ahler-Einstein surface. If the initial compact
surface is symplectic, then the $\lmd$ tangent cone $\g1^\infty$
is stationary.
\end{proposition}
{\it Proof:} Let $V_t^k$ be the varifold defined by $\g1_t^k$. By
the definition of varifolds, we have
$$
V_t^k(\psi)=\int_{\g1_t^k}\psi(x, T{\g1_t^k})d\mu_t^k
$$
for any $\psi\in C_0^0(G^2({\bf R}^4),R)$, where $G^2({\bf R}^4)$
is the Grassmanian bundle of all 2-planes tangent to
$\g1_t^\infty$ in ${\bf R}^4$. For each smooth surface $\g1^k_t$,
the first variation $\d1V_t^k$ of $V_t^k$
(cf. [A], (39.4) in [Si1] and (1.7) in {I2]) is
$$
\d1V_t^k=-\mu_t^k\lfloor \h2_k.
$$
By Proposition \ref{main2.1}, we have that $\d1V_t^k\to 0$ at $t$
as $k\rightarrow\infty$.

By (iii) in Theorem \ref{cmpt}, we have that
$$
-\mu^\infty\lfloor \h2_\infty= \d1 V^\infty =\lim_{k\to\infty} \d1
V_t^k=0.
$$
Therefore $\g1^\infty$ is stationary. \hfill Q.E.D.

\begin{theorem}\label{hol}
Let $M$ be a compact K\"ahler-Einstein surface. If the initial compact surface is
symplectic and $T>0$ is the first blow-up time of the mean
curvature flow, then the $\lmd$ tangent cone $\g1^\infty$ of the
mean curvature flow at $(X_0,T)$ is a holomorphic subvariety of complex dimension one
in some complex structure on ${\bf R}^4$, with
multiplicity more than one in ${\bf R}^4$.
\end{theorem}
{\it Proof:} For a smoothly immersed real surface $\Sigma$ in a
K\"ahler manifold $(M,\omega)$ of complex dimension two, we may
choose a local orthogonal frame $\{e_1,e_2,e_3,e_4\}$ on $M$ along
$\Sigma$ with $e_1,e_2$ tangent to $\Sigma$, so that along
$\Sigma$ the K\"ahler form $\omega$ has the expression (cf. [CT1,
CW]):
$$
\omega=\cos\alpha(e^*_1\wedge e^*_2+e^*_3\wedge e^*_4)
+\sin\alpha(e^*_1\wedge e^*_3-e^*_2\wedge e^*_4),
$$
where $e^*_j$ is the dual of $e_j$ for $j=1,2,3,4$.

Along each surface $\Sigma^k_t$, we have \bs
\omega&=&\cos\alpha_t^k(e^*_1(\g1_t^k)\wedge
e^*_2(\g1_t^k)+e^*_3(\g1_t^k)\wedge e^*_4(\g1_t^k))\\
&& +\sin\alpha_t^k(e^*_1(\g1_t^k)\wedge
e^*_3(\g1_t^k)-e^*_2(\g1_t^k)\wedge e^*_4(\g1_t^k)), \es and
$$
\omega |_{\g1_t^k}=\cos\al_k d\mu_t^k =\cos\alpha_k
e^*_1(\Sigma^k_t)\wedge e^*_2(\Sigma^k_t).
$$

Since $\cos\alpha_k>c>0$ along the flow, the compactness theorem
for Radon measures (cf. Theorem 4.4 in [Si]) implies that the
bounded positive measures $\cos\alpha_k d\mu^k_t$ and $d\mu^k_t$
converge to nonnegative measures $\t1_0d\mu^\infty$ and
$d\mu^\infty$ respectively on $\Sigma^\infty$ in the sense of
measures, for some measurable function $\theta_0$ on
$\Sigma^\infty$ with $0<\theta_0\leq 1$. Here we take
a subsequence if necessary.

In the scaling process on a small neighborhood of the singular
point $X_0$ in $M$, the Riemannian metric $g^k$ tends to the flat
metric on ${\bf R}^4$ as $k\to\infty$, equivalently, the K\"ahler
form $\omega^k$ converges to a self-dual, positive definite,
constant 2-form $\omega_0$ on ${\bf R}^4$ with $\omega_0(0)
=\omega(X_0)$, where $0$ is the origin of ${\bf R}^4$.

By (ii) in Theorem \ref{cmpt}, along $\g1_t^k$, $$ \omega^k
|_{\g1_t^k}\to \omega_0|_{\g1^\infty}=\t1_0\,
e^*_1(\g1^\infty)\wedge e^*_2(\g1^\infty),$$ as measures. Note
that Allard's compactness theorem only provides convergence of
tangent planes to $\g1^k_t$ so the other components in $\omega^k$ may not
converge to those in $\omega_0$. Nevertheless, we have
$$
\theta_0 d\mu^\infty =\theta_0 e^*_1(\Sigma^\infty)\wedge
e^*_2(\Sigma^\infty)= \omega_0 |_{\g1^\infty}.
$$

Next, we shall show that $\t1_0$ is constant ${\cal H}^2$ a.e. on
$\g1^\infty$. To do so, we claim that for any $r>0$,
$\xi_1,\xi_2\in \g1_t^k\cap B_{R/2}(0)$ the following holds
 $$\left|\fr{1}{{\rm
 Vol}(B_r(\xi_1)\cap\g1_t^k)}\int_{B_r(\xi_1)\cap\g1_t^k}\cos\al_k
 d\mu_t^k-\fr{1}{{\rm
 Vol}(B_r(\xi_2)\cap\g1_t^k)}\int_{B_r(\xi_2)\cap\g1_t^k}\cos\al_k
 d\mu_t^k \right|
 $$
 \begin{equation}\label{t1c}
 \leq \fr{C_1(r)}{{\rm
 Vol}(B_r(\xi_1)\cap\g1_t^k)}\cdot \fr{C_2(r)}{{\rm
 Vol}(B_r(\xi_2)\cap\g1_t^k)}\int_{B_R(0)\cap\g1_t^k}\left |\btd
 \cos\al_k
 \right |d\mu_t^k,
 \end{equation}
where $B_r(\xi_i)$, $i=1,2$, are the 4-dimensional balls in $M$.
To prove (\ref{t1c}), let us first recall the isoperimetric
inequality on $\g1_t^k$ (c.f. [HSp] and [MS]): let $B^k_\rho(p)$
be the geodesic ball in $\Sigma^k_t$, with radius $\rho$ and
center $p$, then \bs {\rm Vol}(B^k_\rho(p))&\leq & C\left({\rm
length }(\partial
(B^k_\rho(p)))+\int_{B^k_\rho(p)}|\h2_k|d\mu_t^k\right)^2\\
&\leq & C\left({\rm length }(\partial
(B^k_\rho(p)))+\left(\int_{B^k_\rho(p)}|\h2_k|^2d\mu_t^k
\right)^{1/2}{\rm
Vol}^{1/2}(B^k_\rho(p))\right)^2, \es for any $p\in \g1_t^k$, and
any $\rho >0$, where $C$ does not depend on $k$, $\rho$, and $p$.
By Proposition \ref{main2.1}, we have
$$
\int_{B^k_\rho(p)}|\h2_k|^2d\mu_t^k \to 0~{\rm as}~k\to\infty .
$$
So, for $k$ sufficiently large, we obtain:
$$
{\rm Vol}(B^k_\rho(p))\leq  C\left({\rm length }(\partial
(B^k_\rho(p)))\right)^2.
$$
In particular, for $k$ sufficiently large, the isoperimetric
inequality implies
\begin{equation}\label{iso}
{\rm Vol}(B^k_\rho(p))\geq C\rho^2,
\end{equation}
where $C$ is a positive constant independent of $k,\rho$ and $p$.

Suppose that the diameter of $B_r(\xi)\cap\g1_t^k$ is $d_k(\xi)$.
Then
 \bs
 Cr^2 &\geq & \int_{B_r(\xi)\cap\g1_t^k}d\mu_t^k~~~{\rm by}~(\ref{fm})\\
 &=&\int_0^{d_k(\xi)/2}\int_{\partial B^k_\rho(p)}d\sigma d\rho ~~~{\rm
 for~some}~p\in\g1_t^k\\
 &\geq & c \int_0^{d_k(\xi)/2}{\rm
 Vol}^{1/2}(B^k_\rho(p))d\rho+o(1),~o(1)\to 0~{\rm as}~k\to\infty\\
 &\geq& c\int^{d_k(\xi)/2}_0 C\rho d\rho+o(1)~~{\rm by}~(\ref{iso})\\
 &\geq & cd_k(\xi)^2+o(1).
 \es
We therefore have, for any $\xi$, \be \label{t1c2} d_k(\xi)\leq Cr+o(1) \ee
where the constant $C$ is independent of $\xi$ and $k$.

For any fixed $\eta\in B_r(\xi_2)\cap\g1_t^k$ and any $\xi\in
B_r(\xi_1)\cap\g1_t^k$, we choose a geodesic $l_{\eta\xi}$
connecting $\eta$ and $\xi$, call it a ray from $\eta$ to $\xi$.
Take an open tubular neighborhood $U(l_{\eta\xi})$ of
$l_{\eta\xi}$ in $\g1_t^k$. Within this neighborhood
$U(l_{\eta\xi})$, we call the line in the normal direction of the
ray $l_{\eta\xi}$ the normal line which we denote by
$n(l_{\eta\xi})$. It is clear that \be\label{by parts}
\cos\al_k(\xi)-\cos\al_k(\eta)= \int_{l_{\eta\xi}}\partial_l
\cos\al_k dl \ee where $dl$ is the arc-length element of
$l_{\eta\xi}$. Choose $r$ small enough so that
$B_r(\xi_1)\cap\g1^k_t$ is contained in $U(l_{\eta\xi_1})$.
Keeping $\eta$ fixed and integrating (\ref{by parts}) with respect
to the variable $\xi$, first along the normal direction
$n(l_{\eta\xi_1})$ and then on the ray direction $l_{\eta\xi_1}$,
we have
\begin{eqnarray}\label{average}
\lefteqn{\left|\fr{1}{{\rm
 Vol}(B_r(\xi_1)\cap\g1_t^k)}\int_{B_r(\xi_1)\cap\g1_t^k}\cos\al_k(\xi)
 d\mu_t^k-\cos\al_k(\eta)\right|}\nonumber\\
 &\leq & \fr{1}{{\rm
Vol}(B_r(\xi_1)\cap\g1_t^k)}\int_{0}^{d_k(\xi_1)}\int_{n(l_{\eta\xi_1})}\int_{l_{\eta\xi}}
\left |\btd \cos\al_k
 \right |dl d n(\xi)d\rho\nonumber\\
 &\leq & \fr{1}{{\rm
Vol}(B_r(\xi_1)\cap\g1_t^k)}\int_{0}^{d_k(\xi_1)}\int_{B_R(0)}
\left |\btd \cos\al_k
 \right |d\mu_t^kd\rho\nonumber\\
&\leq&\fr{Cr}{{\rm Vol}(B_r(\xi_1)\cap\g1_t^k)}\int_{B_R(0)}\left
|\btd \cos\al_k \right |d\mu_t^k,
\end{eqnarray}
here in the last step we have used (\ref{t1c2}). From
(\ref{average}), integrating with respect to $\eta$ in
$B_r(\xi_2)\cap\g1_t^k$ and dividing by ${\rm
Vol}(B_r(\xi_2)\cap\g1^k_t)$, we get the desired inequality
(\ref{t1c}).

For $i=1,2$ H\"older's inequality and (\ref{fm}) lead to
$$
\int_{B_r(\xi_i)\cap\g1^k_t}\left|\btd\cos\alpha_k\right|d\mu^k_t\leq
Cr\left(\int_{B_r(\xi_i)\cap\g1^k_t}\left|\btd\cos\alpha_k\right|^2d\mu^k_t\right)^{1/2}.
$$
The triangle inequality implies $B^k_r(\xi_i)\subset
B_r(\xi_i)\cap\g1^k_t$ for $i=1,2$; therefore by (\ref{iso})
$$
{\rm Vol}(B_r(\xi_i)\cap \g1^k_t)\geq {\rm Vol}(B^k_r(\xi_i))\geq
Cr^2.
$$
Now first letting $k\to\infty$ in (\ref{t1c}) and using that the
right hand side of (\ref{t1c}) tends to 0 by Proposition
\ref{main2.1}, and then letting $r\to 0$, we conclude that $\t1$ is
constant ${\cal H}^2$ a.e. on $\g1^\infty$.

We choose a new complex structure on ${\bf R}^4$ so that with
respect to the corresponding K\"ahler form $\omega_1$ we have
$$ \omega_1 |_{\g1^\infty}=\t1_1d\mu^\infty ~{\rm
and}~\t1_1=1~{\rm at~one~point}~x\in\g1^\infty .$$ It is clear
that, in the new complex structure on ${\bf R}^4$, Proposition
\ref{main2.1} remains true. Therefore, by the same argument as above for $\theta$,
we conclude that $\theta_1\equiv 1$ and accordingly
 \be\label{hol3} \omega_1 |_{\g1^\infty}\equiv d\mu^\infty
.\ee

In other words, we have shown that $\g1^\infty$ is calibrated by the
closed 2-form $\omega_1$. Applying Theorem 2.1 in [HS], we see that
$\g1^\infty$ is a holomorphic subvariety, in the chosen complex structure,
 of complex dimension one. More
precisely, Harvey-Shiffman's theorem says that, if $T$ is a
locally rectifiable $(k,k)$ type current and the $(2k+1)$-dim
measure of ${\rm supp}(T)=0$, then $T$ is a holomorphic $k$-chain,
i.e., away from the support of its boundary $T$ is a locally
finite sum: $T=\sum_j n_j[V_j]$, where $n_j$ are positive integers
and $V_j$ are pure $1$-dimensional complex subvarieties (cf. [HS],
[HL]) . We proved that $\g1^\infty$ is ${\cal H}^2$-rectifiable
(see Proposition \ref{rectif}) and it is stationary (see Lemma
\ref{stat}), which implies that the $3$-dimensional measure of
${\rm supp}(\g1^\infty)=0$ (cf. Remark. after theorem 2.1 in
[HS]). By (\ref{hol3}), we can see that $\g1^\infty$ is a bidegree $(1,1)$
current, hence a holomorphic curve according to
Harvey-Shiffman's theorem.

We are left to show that the holomorphic curve has multiplicity more than
1. Otherwise, we would have
$$
\lim_{\rho\to 0}\fr{1}{\pi\rho^2}\mu^\infty(B_\rho(0))=1.
$$
It then follows from (\ref{smon}) that for any $\e1>0$, there are
$\d1>0$ and $k_0>0$ such that for any $0<\rho<2\d1$ and $k>k_0$,
 \be\label{rec11} \rho^{-2}\mu_{-r^2}^k(B_\rho(\xi))< \pi
(1+\e1)
\ee
for any fixed $r>0$. The choice of $r$ will be based on
the following observation
\bs
\Phi(F_k,0,0,0-r^2)&\leq &\fr{1}{4\pi r^2}\int_0^{\d1}\int_{\part
B_\rho(0)\cap\g1_{0-r^2}^k}e^{-\fr{\rho^2}{4r^2}}d\mu_{0-r^2}^k\\
&&+\fr{1}{4\pi r^2}\int_{\d1}^{\infty}\int_{\part
B_\rho(0)\cap\g1_{0-r^2}^k}e^{-\fr{\rho^2}{4r^2}}d\mu_{0-r^2}^k\\
&\leq &\fr{1}{4\pi r^2}\int_0^{\d1} e^{-\fr{\rho^2}{4r^2}}\int_{
\partial B_\rho(0)\cap\g1_{0-r^2}^k}d\mu_{0-r^2}^k\\
&&+\fr{1}{4\pi r^2}\int_{\d1}^\infty e^{-\fr{\rho^2}{4r^2}}\int_{
\partial B_\rho(0)\cap\g1_{0-r^2}^k}d\mu_{0-r^2}^k\\
&\leq &\fr{1}{4\pi r^2}\int_0^{\d1}
e^{-\fr{\rho^2}{4r^2}}\fr{\rho}{2r^2}{\rm Vol}(
B_\rho(0)\cap\g1_{0-r^2}^k)d\rho\\
&&+\fr{1}{4\pi r^2}\int_{\d1}^\infty
e^{-\fr{\rho^2}{4r^2}}\fr{\rho}{2r^2}{\rm Vol}(
B_\rho(0)\cap\g1_{0-r^2}^k)d\rho\\
&\leq &\fr{1}{4\pi r^2}\int_0^\d1
e^{-\fr{\rho^2}{4r^2}}\fr{\rho^3}{2r^2}d\rho +\e1 + o(r)~~{\rm
by}~~(\ref{rec11})~~{\rm and}~~(\ref{fm})\\
&\leq & 1 +o(r).
\es
Choosing $r>0$ sufficiently small, we therefore have
$$
\Phi (F,X_0,T,T-\lmd_k^{-2}r^2)=\Phi(F_k,0,0,0-r^2)\leq 1+\e1.
$$
Now by White's local regularity theorem ([Wh1] Theorem 3.1
and Theorem 4.1, also see [E2]), $(X_0,T)$ could not be
a singular point of the mean curvature flow. This is a contradiction.
Now the proof of Theorem \ref{hol} is complete. \hfill Q.E.D.

\section{Flatness of the $\lmd$ tangent cones}

In this section, we prove that the $\lmd$ tangent cones are flat.
\begin{theorem}\label{flat}
Let $M$ be a compact K\"ahler-Einstein surface. If the initial compact
surface is symplectic and $T>0$ is the first blow-up time of the
mean curvature flow, then the $\lmd$ tangent cone $\g1^\infty$ of
the mean curvature flow at $(X_0,T)$ consists of a finite union of
more than one 2-planes which are complex in a complex structure on
${\bf R}^4$ .
\end{theorem}
{\it Proof:} Without loss of any generality, we may assume
$0\in\g1^\infty$ where $0$ is the origin of ${\bf R}^4$.
In fact, if not, $\g1^\infty$ would
move to infinity, then we would have
 $$
\Phi (F,X_0,T,T-\lmd_k^{-2}r^2)=\Phi(F_k,0,0,0-r^2)\to 0~{\rm
as}~k\to\infty .
$$
But White's regularity theorem then implies that $(X_0,T)$ is a
regular point. This is impossible.

There is a sequence of points $X_k\in\g1_t^k$ satisfying $X_k\to
0$ as $k\to\infty$. By Proposition \ref{main2.1}, for any $s_1$ and $s_2$ with
$-\infty<s_1<s_2<0$ and any $R>0$, we have
$$
\int_{s_1}^{s_2}\int_{\g1_t^k\cap
B_R(0)}\left|F_k^{\perp}\right|^2d\mu_t^k dt\to 0~~{\rm
as}~~k\to\infty.
$$
Thus, by (\ref{fm})
 \bs
\lefteqn{\lim_{k\to\infty}\int_{s_1}^{s_2}\int_{\g1_t^k\cap
B_R(0)}\left|(F_k-X_k)^{\perp}\right|^2d\mu_t^k dt}\\
&\leq & 2\lim_{k\to\infty}\int_{s_1}^{s_2}\int_{\g1_t^k\cap
B_R(0)}\left|F_k^{\perp}\right|^2d\mu_t^k dt+ C(s_2-s_1)R^2
\lim_{k\to\infty}|X_k|^2\\
&=& 0.
 \es

Let us denote the tangent spaces of $\g1^k_t$ at the point $F_k(x,t)$
and of $\g1^\infty$ at the point $F^\infty(x,t)$ by $T\g1_t^k$ and $T\g1^\infty$
respectively. It is clear that
$$
(F_k-X_k)^{\perp}={\rm dist}~(X_k, T{\g1_t^k}),
$$
and
$$
(F_\infty)^{\perp}={\rm dist}~(0, T{\g1^\infty}).
$$
By Allard's compactness theorem, i.e. Theorem \ref{cmpt} (ii), we
have \bs \int_{s_1}^{s_2}\int_{\g1^\infty\cap
B_R(0)}\left|(F_\infty)^{\perp}\right|^2d\mu^\infty dt &=&
\int_{s_1}^{s_2}\int_{\g1^\infty\cap B_R(0)}\left|{\rm dist}~(0,
T{\g1^\infty})\right|^2d\mu^\infty dt\\
&=& \lim_{k\to\infty}\int_{s_1}^{s_2}\int_{\g1_t^k\cap
B_R(0)}\left|{\rm dist}~(X_k, T{\g1_t^k})\right|^2d\mu_t^k dt\\
&=& \lim_{k\to\infty}\int_{s_1}^{s_2}\int_{\g1_t^k\cap
B_R(0)}\left|(F_k-X_k)^{\perp}\right|^2d\mu_t^k dt\\
&=&0.\es

By Theorem \ref{hol}, we know that each sheet of $\g1^\infty$ is
smooth outside a discrete set of points $\cal S$. So outside $\cal
S$, we have
$$
\la F_\infty,v_\al\ra=0.
$$
Note that the above inner product is taken in ${\bf R}^4$, and
differentiating in ${\bf R}^4$ then yields
$$
0=\la \p_i F_\infty,v_\al\ra+\la F_\infty ,\p_i v_\al\ra =\la
F_\infty ,\p_i v_\al\ra.
$$
Because $\p_i F_\infty$ is tangential to $\Sigma^\infty$,
by Weingarten's equation we observe
$$
(h_\infty)_{ij}^\al\la\f2_\infty,e_j\ra=0~~{\rm for}~~{\rm
all}~~\al,~i=1,2.
$$
So for $\al=1,2$, we have
$$
\det ((h_\infty)_{ij}^\al)=0.
$$
Since $\h2 =0$, for $\al=1,2$ we also have
$$
tr ((h_\infty)_{ij}^\al)=0.
$$
It then follows immediately that the symmetric matrix $((h_\infty)_{ij}^\al)$
is in fact the zero matrix, for all $i,j,\al=1,2$, which obviously yields
$|\a2_\infty|\equiv 0$.

Since the second fundamental form $\a2_\infty$ of $\Sigma^\infty$
is identically zero on the smooth locus of $\Sigma^\infty$ whose
multiplicity is finite but bigger than 1,
$\Sigma^\infty$ is a finite union of more than one 2-planes.
Moreover, if any two of these 2-planes meet at two distinct points
they would intersect along a line containing these two points; but
this contradicts to that $\Sigma^\infty$ is a holomorphic
subvariety.

 This completes the proof of Theorem \ref{flat}. \hfill{Q.E.D.}

\section{Tangent cones from time dependent scaling}
In this section, we consider the tangent cones which arise from
the rescaled surface $\tilde{\g1}_s$ defined by
\begin{equation}
\tilde{F}(\cdot ,s)=\fr{1}{\sqrt{2(T-t)}}F(\cdot ,t),
\end{equation}
where $s=-\fr{1}{2}\log (T-t)$, $c_0\leq s<\infty $. Here we choose the
coordinates so that $X_0=0$. Rescaling of
this type was used by Huisken [H2] to distinguish Type I and Type
II singularities for mean curvature flows. Denote the rescaled surface by
$\tilde{\g1}_s$. From the evolution equation of $F$ we derive
the flow equation for $\tilde{F}$
\begin{equation}\label{meaneqn2}
\fr{\p}{\p s}\tilde{F}(x,s)=\tilde{\h2}(x,s)+\tilde{\f2}(x,s).
\end{equation}
It is clear that
\begin{eqnarray*}
\cos\tilde{\al}(x,s)&=&\cos\al(x,s),\\
|\tilde{\a2}|^2(x,s)&=&2(T-t)|\a2 |^2(x,t).
\end{eqnarray*}
Recall that we set, in section 2, $v(x,t)=e^{R_0t}\cos\al(x,t)$ where
$R_0=\max\{0,-R\}$ and $R$ is the scalar curvature of $M$. The function $v(x,t)$
satisfies
$$
\left(\fr{\p}{\p t}-\btu \right)v(x,t) \geq
\left|\overline{\nabla}J_{\g1_t}\right|^2v(x,t).
$$
The corresponding version of this evolution inequality for the scaled flow is
in the following lemma.
\begin{lemma}\label{tkae}
Assume that $M$ is a K\"ahler-Einstein surface with scalar
curvature $R$ and $\g1_t$ evolves by a mean curvature flow in $M$ with the initial surface
being compact and symplectic. Let
$\tilde{v}(x,s)=e^{R_0(T-e^{-2s})}\cos\tilde{\al}(x,s)$ where $R_0=\max\{0,-R\}$. Then
\be\label{tkae1}
\left(\fr{\p}{\p s}-\tilde{\btu} \right)\tilde{v}(x,s) \geq |\overline{\nabla}
 J_{\tilde{\g1}_s}|^2\tilde{v}(x,s).\ee
\end{lemma}
{\it Proof:} One can check directly that
$$
\left(\fr{\p}{\p s}-\tilde{\btu}
\right)\cos\tilde{\al}(x,s)=2(T-t)\left(\fr{\p}{\p t}-\btu
\right)\cos\al(x,t).
$$
It follows that \bs \left(\fr{\p}{\p s}-\tilde{\btu}
\right)\tilde{v}(x,s)&=&
2(T-t)(\fr{\p}{\p t}-\btu )v(x,t)\\
&\geq&2(T-t)\left|\overline{\nabla}J_{\Sigma_t}\right|^2 v(x,t)\\
&=& |\overline{\nabla}
 J_{\tilde{\g1}_s}|^2\tilde{v}(x,s).
\es
 This proves the lemma. \hfill{Q.E.D.}

Next, we shall derive the corresponding weighted monotonicity formula
for the scaled flow. By (\ref{tkae1}), we have
$$
\left(\fr{\p}{\p s}-\tilde{\btu} \right)\fr{1}{\tilde{v}}\leq
-|\overline{\nabla}J_{\tilde{\g1_t}}|^2\fr{1}{\tilde{v}}-\fr{2}{\tilde{v}^3}
|\tilde{\btd}{\tilde{v}}|^2.
$$
Let
$$
\tilde{\rho}(X)={\rm exp}\left(-\fr{1}{2}|X|^2\right),
$$
$$
\Psi (s)=\int_{\tilde{\g1}_s}\fr{1}{\tilde{v}}\phi\tilde{\rho}
(\tilde{\f2})d\tilde{\mu}_s.
$$
\begin{lemma} There are positive constants $c_1$ and $c_2$ which
depend on $M$, $\f2_0$ and $r$ which is the constant in the
definition of $\phi$, so that the following monotonicity formula
holds \ba \label{mon3} \fr{\p}{\p s}{\rm exp} (c_1 e^{-s})\Psi (s)
&\leq&-{\rm exp} (c_1 e^{-s}) \left(
\int_{\tilde{\g1}_s}\fr{1}{\tilde{v}}\phi\tilde{\rho}
(\tilde{\f2})
\left|\tilde{\h2}+ \tilde{F}^{\perp}\right|^2d\tilde{\mu}_s \right. \no\\
&&\left.+ \int_{\tilde{\g1}_s}\fr{1}{\tilde{v}}\phi\tilde{\rho}
(\tilde{\f2})\left|\overline{\nabla}J_{\tilde{\g1}_s}\right|^2d\tilde{\mu}_s
+\int_{\tilde{\g1}_s}\fr{2}{\tilde{v}^3}\left|\tilde{\btd}
\tilde{v}\right|^2\phi\tilde{\rho} (\tilde{\f2})d\tilde{\mu}_s
\right)\no\\
&&+c_2e^{-2s}. \ea
\end{lemma}
{\it Proof:} Note that
\begin{eqnarray*}
&&\tilde{F}(x,s)=\fr{F(x,t)}{\sqrt{2(T-t)}},\\
&&\tilde{\h2}(x,s)=\sqrt{2(T-t)}\h2(x,t),\\
&&|\overline{\nabla}J_{\tilde{\g1}_s}|^2(x,s)
=2(T-t)|\overline{\nabla}J_{\g1_t}|^2(x,t),\\
&&|\tilde{\btd}\tilde{v}|^2(x,s)=2(T-t)|\btd v|^2(x,t).
\end{eqnarray*}
By the chain rule
$$\fr{\p}{\p s}=2e^{-2s}\fr{\p}{\p t}$$
and the monotonicity inequality (\ref{mon2}) for the unscaled
surface, we obtain the desired inequality. \hfill{Q.E.D.}

\begin{lemma}
\label{main7.1} Let $M$ be a compact K\"ahler-Einstein surface. If
the initial compact surface is symplectic, then there is a sequence
$s_k\to\infty$ such that, for any $R>0$, \be\label{m7.1}
\int_{\tilde{\g1}_{s_k}\cap B_R(0)}
|\bar{\btd}J_{\tilde{\g1}}|^2d\tilde{\mu}_{s_k} \to 0~~{\rm
as}~~k\to\infty, \ee \be\label{m7.2} \int_{\tilde{\g1}_{s_k}\cap
B_R(0)} |\btd \cos \tilde{\al}|^2d\tilde{\mu}_{s_k}\to 0~~{\rm
as}~~k\to\infty, \ee \be\label{m7.3} \int_{\tilde{\g1}_{s_k}\cap
B_R(0)}|\tilde{\h2}|^2d\tilde{\mu}_{s_k}\to 0~~{\rm
as}~~k\to\infty, \ee and \be\label{m7.4}
\int_{\tilde{\g1}_{s_k}\cap
B_R(0)}|\tilde{F}^{\perp}|^2d\tilde{\mu}_{s_k}\to 0~~{\rm
as}~~k\to\infty. \ee
\end{lemma}
{\it Proof:} By (\ref{mon3}), we have \bs  \infty & > &
\int_{s_0}^\infty\int_{\tilde{\g1}_s}\fr{1}{\tilde{v}}\phi\tilde{\rho}
(\tilde{\f2})
\left|\tilde{\h2}+ \tilde{F}^{\perp}\right|^2d\tilde{\mu}_s ds \\
&&+
\int_{s_0}^\infty\left(\int_{\tilde{\g1}_s}\fr{1}{\tilde{v}}\phi\tilde{\rho}
(\tilde{\f2})|\overline{\nabla}J_{\tilde{\g1}_s}|^2d\tilde{\mu}_s
+\int_{\tilde{\mu}_s}\fr{2}{\tilde{v}^3}|\tilde{\btd}
\tilde{v}|^2\phi\tilde{\rho} (\tilde{\f2})d\tilde{\mu}_s
\right)ds. \es
Hence there is a sequence $s_k\to\infty$, such that as $k\to\infty$
$$
\int_{\tilde{\g1}_{s_k}}\fr{1}{\tilde{v}}\phi\tilde{\rho}
(\tilde{\f2})|\overline{\nabla}J_{\tilde{\g1}_{s_k}}|^2d\tilde{\mu}_{s_k}\to
0
$$
and
$$
\int_{\tilde{\g1}_{s_k}}\fr{2}{\tilde{v}^3}|\tilde{\btd}
\tilde{v}|^2\phi\tilde{\rho} (\tilde{\f2})d\tilde{\mu}_{s_k}\to 0,
$$
which yields (\ref{m7.1}) and (\ref{m7.2}), and \be\label{m7.5}
\int_{\tilde{\g1}_{s_k}}\fr{1}{\tilde{v}}\phi\tilde{\rho}
(\tilde{\f2}) \left|\tilde{\h2}+
\tilde{F}^{\perp}\right|^2d\tilde{\mu}_{s_k}\to 0. \ee By
(\ref{djh}) and (\ref{m7.1}) we see that (\ref{m7.3}) holds.
By (\ref{m7.3}) and (\ref{m7.5}) we get (\ref{m7.4}). This proves
the proposition. \hfill Q.E.D.

The proof of the following lemma is essentially the same as the
one for Proposition \ref{rectif}, except there are two parameters
$\lmd,t$ for the $\lmd$ tagent cones but only one parameter $t$
for the time dependent tangent cones. For the sake of
completeness, we shall provide a proof.

\begin{lemma}\label{texit}
There is a subsequence of $s_k$, which we also denote by $s_k$,
such that $(\tilde{\g1}_{s_k},d\tilde{\mu}_{s_k})\to
(\tilde{\g1}_\infty,d\tilde{\mu}_\infty)$ in the sense of measures.
And $(\tilde{\g1}_\infty,d\tilde{\mu}_\infty)$ is ${\mathcal H}^2$-rectifiable.
\end{lemma}
{\it Proof:} To show the subconvergence, it suffices to show that,
for any $R>0$, \be\label{tfm}
\tilde{\mu}_{s_k}(\tilde{\g1}_{s_k}\cap B_R(0))\leq CR^2,
 \ee
where $B_R(0)$ is a metric ball in ${\bf R}^4$, $C>0$ is
independent of $k$. Direct calculation leads to
\bs
\tilde{\mu}_{s_k}(\tilde{\g1}_{s_k}\cap B_R(0))
&=&(2(T-t))^{-1}\int_{\g1_{T-e^{2s_k}}\cap B_{\sqrt{2(T-t)}R}(0)}d\mu_t\\
&=& R^2\left(\sqrt{2}e^{-s}R\right)^{-2}\int_{\g1_{T-e^{2s_k}}\cap
B_{\sqrt{2}
e^{-s}R}(0)}d\mu_t\\
&\leq & CR^2\Phi\left(0,T+(\sqrt{2}e^{-s}R)^{2}-e^{2s_k},
T-e^{2s_k}\right) \es By the monotonicity inequality (\ref{mon1}),
we have \bs \tilde{\mu}_{s_k}(\tilde{\g1}_{s_k}\cap B_R(0))&\leq
&CR^2\left(\Phi(0,
T+(\sqrt{2}e^{-s}R)^{2}-e^{2s_k}, T/2)+C\right)\\
&\leq &C\fr{R^2}{T}\left(\mu_{T/2}(\g1_{T/2})+C\right). \es Since
$$
\fr{\p}{\p t}\mu_t(\g1_t)=-\int_{\g1_t}|\h2|^2d\mu_t,
$$
we have
$$
\tilde{\mu}_{s_k}(\tilde{\g1}_{s_k}\cap B_R(0))\leq CR^2.
$$

We now prove that $(\tilde{\g1}_\infty,d\tilde{\mu}_\infty)$ is
${\mathcal H}^2$-rectifiable. For any $\xi\in\tilde{\g1}_\infty$,
choose $\xi_k\in\tilde{\g1}_{s_k}$ with $\xi_k\to\xi$ as
$k\to\infty$. By the monotonicity identity (17.4) in [Si1], we
have \ba\label{tsmon}
\sigma^{-2}\tilde{\mu}_{s_k}(B_\sigma(\xi_k))&=&
\rho^{-2}\tilde{\mu}_{s_k}(B_\rho(\xi_k))
-\int_{B_\rho(\xi_k)\setminus B_\sigma(\xi_k)}\fr{|D^\perp
r|^2}{r^2} d\tilde{\mu}_{s_k}
\no\\
&& -\fr{1}{2}\int_{B_\rho(\xi_k)}(x-\xi_k)\cdot \tilde{{\bf H}_k}
\left(\fr{1}{r_\sigma^2}-\fr{1}{\rho^2}\right)d\tilde{\mu}_{s_k},
\ea for all $0<\sigma \leq \rho$, where
$\tilde{\mu}_{s_k}(B_\sigma(\xi_k))$ is the area of
$\tilde{\g1}_{s_k}\cap B_\sigma(\xi_k)$, $r_\sigma =\max
\{r,\sigma\}$ and $D^\perp r$ denotes the orthogonal projection of
$Dr$ (which is a vector of length 1) onto
$(T_{\xi_k}\tilde{\g1}_{s_k})^\perp$. Letting $k\to\infty$, by
Lemma \ref{main7.1}, we have
$$
\sigma^{-2}\tilde{\mu}_\infty(B_\sigma(\xi))\leq
\rho^{-2}\tilde{\mu}_\infty(B_\rho(\xi)),
$$
for all $0<\sigma \leq \rho$. Therefore, $\lim_{\rho\to
0}\rho^{-2}\tilde{\mu}_\infty(B_\rho(\xi))$ exists. By (\ref{tfm})
we know that $\lim_{\rho\to 0}\rho^{-2}\tilde{\mu}_\infty(B_\rho(\xi))$
is finite; and we now show that it has a positive lower bound.

By the isoperimetric inequality on $\tilde{\g1}_{s_k}$ (c.f. [HSp]
and [MS]), we have \bs {\rm Vol}(B^k_\rho(\xi_k))&\leq &
C\left({\rm length }(\partial
(B^k_\rho(\xi_k)))+\int_{B^k_\rho(\xi_k)}|\tilde{\h2}|d\tilde{\mu}_{s_k} \right)^2\\
&\leq & C\left({\rm length }(\partial
(B^k_\rho(\xi_k)))+\left(\int_{B^k_\rho(\xi_k)}|\tilde{\h2}|^2d\tilde{\mu}_{s_k}
\right)^{1/2}{\rm Vol}^{1/2}(B^k_\rho(\xi_k))\right)^2, \es for
 any $\rho >0$, where
$B^k_\rho(\xi_k)$ is the geodesic ball in $\tilde{\g1}_{s_k}$,
with radius $\rho$ and center $\xi_k$, $C$ does not depend on $k$,
$\rho$. By Lemma \ref{main7.1}, we have
$$
\int_{B^k_\rho(\xi_k)}|\tilde{\h2}|^2d\tilde{\mu}_{s_k} \to 0~{\rm
as}~k\to\infty .
$$
Hence, for $k$ sufficiently large, we have
\begin{equation}\label{tiso}
{\rm Vol}(B^k_\rho(\xi_k))\geq C\rho^2,
\end{equation}
where $C$ is a positive constant independent of $k$, $\rho$.

The triangle inequality implies $B^k_r(\xi_k)\subset
B_r(\xi_k)\cap\tilde{\g1}_{s_k}$ for $k=1,2,\cdots$, therefore by
(\ref{tiso})
$$
{\rm Vol}(B_r(\xi_k)\cap \tilde{\g1}_{s_k})\geq {\rm
Vol}(B^k_r(\xi_k))\geq Cr^2.
$$

It concludes that $\lim_{\rho\to
0}\rho^{-2}\tilde{\mu}_\infty(B_\rho(\xi))$ exists and for
${\mathcal H}^2$ almost all $\xi\in \tilde{\g1}_\infty$,
\begin{equation}\label{tden}
0<C\leq\lim_{\rho\to
0}\rho^{-2}\tilde{\mu}_\infty(B_\rho(\xi))<\infty.
\end{equation}

By Priess's theorem in [P] we can see from (\ref{tden}) that,
$(\tilde{\g1}_\infty,d\tilde{\mu}_\infty)$ is ${\mathcal H}^2$-rectifiable.
This proves the lemma. \hfill Q.E.D.

\begin{definition}
\label{bbtt} {\em Let $(X_0,T)$ be a singular point of the mean
curvature flow of a closed symplectic surface in a compact
K\"ahler-Einstein surface $M$. We call
$(\tilde{\g1}_\infty,d\tilde{\mu}_\infty)$ obtained in the last
lemma {\it a tangent cone  of the mean curvature flow $\g1_t$ at
$(X_0,T)$ in the time dependent scaling.}}
\end{definition}

With the lemmas established in this section, we can derive, by using arguments
completely similar to those for the $\lmd$ tangent cones in
section 4 and section 5, holomorphicity and
flatness of the tangent cones coming from time dependent scaling.

\begin{theorem}\label{tflat}
Let $M$ be a compact K\"ahler-Einstein surface. If the initial compact
surface is symplectic and $T>0$ is the first blow-up time of the
mean curvature flow, then the tangent cone $\tilde{\g1}_\infty$ of
the mean curvature flow at $(X_0,T)$ coming from time dependent
scaling consists of a finite union of more than one 2-planes in
${\bf R}^4$ which are complex in a complex structure on ${\bf
R}^4$.
\end{theorem}

Finally, we give two remarks.
\begin{remark} {\em Let $M$ be a compact K\"ahler-Einstein surface.
Assume that the initial surface is symplectic. If $(X_0,T)$ is the
Type I singularity for the mean curvature flow, then
$$
0<c\leq (T-t)|\a2|^2\leq C<\infty ,
$$
and consequently, $\tilde{\g1}_{s_k}$ converges strongly to
$\tilde{\g1}_\infty$ with
$$
0<c\leq |\tilde{\a2}_\infty |^2\leq C<\infty .
$$
However, by Theorem \ref{tflat}, we have $\tilde{\a2}_\infty\equiv
0$. This contradiction shows that there is no Type I singularity
for mean curvature flow of symplectic surface in K-E surfaces.
This result was proved in [CL] and [W1].}
\end{remark}

\begin{remark} \label{morgan}{\em
When one considers the singularity of holomorphic curves, using
the monotonicity identities (\ref{smon}) and (1.2) in [Si2], by an
argument similar to the one used in the present paper, one can
show that the bubbles are two-dimensional planes which are complex
under some complex structure on ${\bf R}^4$ (also see [M]). In
fact, in [M], Morgan obtained more general results, in particular,
he proved that any tangent cone to a two-dimensional oriented area
minimizing surface in ${\bf R}^4$ consists of such planes.}
\end{remark}

\vspace{.2in}
\begin{center}
{\large\bf REFERENCES}
\end{center}
\footnotesize
\begin{description}
\item[{[A]}] {W. Allard, First variation of a varifold, Annals
of Math. {\bf 95} (1972), 419-491.}
\item[{[Ar]}] {C. Arezzo, Minimal surfaces and deformations of holomorphic curves
in K\"ahler-Einstein manifolds, Ann. Scuola Norm. Sup. Pisa Cl. Sci. (4) 29 (2000), no. 2, 473--481.}
\item[{[B]}] {K. Brakke, The motion of a surface by its mean curvature,
Princeton Univ. Press, 1978.}
\item[{[CL]}] {J. Chen and J. Li, Mean curvature flow of surface in 4-manifolds,
Adv. Math. {\bf 163} (2001), 287-309.}
\item[{[CLT]}] {J. Chen, J. Li and G. Tian, Two-dimensional graphs
moving by mean curvature flow, Acta Math. Sinica, English Series
Vol.{\bf 18}, No.2 (2002), 209-224.}
\item[{[CT1]}] {J. Chen and G. Tian, Minimal surfaces in
Riemannian 4-manifolds, GAFA, {\bf 7} (1997), 873-916.}
\item[{[CT2]}] {J. Chen and G. Tian, Moving symplectic curves in
K\"ahler-Einstein surfaces, Acta Math. Sinica, English Series,
{\bf 16} (4), (2000), 541-548.}
\item[{[CW]}] {S.S. Chern and R. Wolfson, Minimal surfaces by
moving frames, Ann. of Math. {\bf 105} (1983), 59-83.}
\item[{[E1]}] {K. Ecker, On regularity for mean curvature flow of hypersurfaces,
Calc. Var. 3(1995), 107-126.}
\item[{[E2]}] {K. Ecker, A local monotonicity formula for mean curvature flow, Ann. Math.,
154(2001), 503-525.}
\item[{[E3]}] {K. Ecker, Lectures on regularity for mean curvature
flow, preprint (2002).}
\item[{[HL]}] {R. Harvey and H.B. Lawson, Calibrated geometries,
Acta Math. {\bf 148} (1982), 47-157.}
\item[{[HS]}] {R. Harvey and B. Shiffman, A characterization of
holomorphic chains, Acta Math. {\bf 148} (1982), 47--157.}
\item[{[HSp]}] {D. Hoffman and J. Spruck, Sobolev and isoperimetric
inequalities for Riemannian submanifolds, Comm. Pure Appl. Math.
{\bf 27} (1974), 715-727.}
\item[{[H1]}] {G. Huisken, Asymptotic behavior for singularities
of the mean curvature flow, J. Diff. Geom. {\bf 31} (1990),
285-299.}
\item[{[H2]}] {G. Huisken, Flow by mean curvature of convex surfaces into
spheres, J. Diff. Geom. {\bf 20} (1984), 237-266.}
\item[{[H3]}] {G. Huisken, Contracting convex hypersurfaces in Riemannian
manifolds by their mean curvature, Invent. math. {\bf 84} (1986),
463-480.}
\item[{[HS1]}] {G. Huisken and C. Sinestrari, Convexity
estimates for mean curvature flow and singularities of mean convex
surfaces. Acta Math. {\bf 183} (1999), no. 1, 45-70.}
\item[{[HS2]}] {G. Huisken and C. Sinestrari, Mean curvature flow
singularities for mean convex surfaces. Calc. Var. Partial
Differential Equations {\bf 8} (1999), no. 1, 1-14.}
\item[{[I1]}] {T. Ilmanen, Singularity of mean curvature flow of surfaces,
preprint.}
\item[{[I2]}] {T. Ilmanen, Elliptic Regularization and Partial Regularity
for Motion by Mean Curvature, Memoirs of the Amer. Math. Soc.,
520, 1994.}
\item[{[M]}] {F. Morgan, On the singular structure of
two-dimensional area minimizing surfaces in ${\bf R}^n$, Math.
Ann. 261(1982), 101-110.}
\item[{[MS]}] {J.H. Michael and L. Simon, Sobolev and mean-valued
inequalities on generalized submanifolds of $R^n$, Comm. Pure
Appl. Math. {\bf 26} (1973), 361-379.}
\item[{[MW]}] {M. Micaleff and B. White, The structure of branch
points in minimal surfaces and in pseudo-holomorphic curves, Ann.
of Math., {\bf 141} (1995), 35-85.}
\item[{[P]}] {D. Priess, Geometry of measures in $R^n$; Distribution,
rectifiability, and densities; Ann. of Math., {\bf 125} (1987),
537-643.}
\item[{[Sh]}] {V. Shevchishin, Pseudo-holomorphic curves and the
symplectic isotopy problem, preprint.}
\item[{[Si1]}] {L. Simon, Lectures on Geometric Measure Theory,
Proc. Center Math. Anal. 3 (1983), Australian National Univ.
Press.}
\item[{[Si2]}] {L. Simon, Existence of surfaces minimizing the
Willmore functional, Comm. Anal. Geom. {\bf 1} (1993), 281-326.}
\item[{[Sm1]}] {K. Smoczyk, Der Lagrangesche mittlere
Kruemmungsfluss. Univ. Leipzig (Habil.-Schr.), 102 S. 2000.}
\item[{[Sm2]}] {K. Smoczyk, Harnack inequality for the Lagrangian
mean curvature flow. Calc. Var. PDE, {\bf 8} (1999), 247-258.}
\item[{[Sm3]}] {K. Smoczyk, Angle theorems for the Lagrangian mean
curvature flow, preprint 2001.}
\item[{[Sp]}] {M. Spivak, A Comprehensive Introduction to
Differential Geometry,Volume 4, Second Edition, Publish or Perish,
Inc. Berkeley, 1979.}
\item[{[ST]}] {B. Siebert and G. Tian, Holomorphy of genus two
 Lefschetz fibration, preprint.}
\item[{[T]}] {G. Tian, Symplectic isotopy in four dimension,
First International Congress of Chinese Mathematicians (Beijing,
1998), 143-147, AMS/IP Stud. Adv. Math., 20, Amer. Math. Soc.,
Providence, RI, 2001. }
\item[{[Wa1]}] {M.-T. Wang, Mean curvature flow of surfaces in
Einstein four manifolds, J. Diff. Geom. {\bf 57} (2001), 301-338.}
\item[{[Wa2]}] {M.-T. Wang, Long-time existence and convergence
of graphic mean curvature flow in arbitrary
codimension, preprint.}
\item[{[Wa3]}] {M.-T. Wang, Deforming area preserving diffeomorphism
of surfaces by mean curvature flow,
 Math. Res. Lett. {\bf 8}, no. 5-6 (2001), 651-661.}
\item[{[Wh1]}] {B. White, A local regularity theorem for classical
mean curvature flow, preprint (2000).}
\item[{[Wh2]}] {B. White, The size of the singular set in mean curvature flow
of mean-convex sets, J. Amer. Math. Soc. {\bf 13}, no.3 (2000),
665--695.}
\end{description}

\end{document}